\documentclass[11pt,reqno]{amsart}
\usepackage[utf8]{inputenc}
\usepackage[T1]{fontenc}
\usepackage{lmodern}
\usepackage{amsmath,amssymb,amsthm,mathtools,mathrsfs}
\usepackage{microtype}
\usepackage{enumitem}
\usepackage[hidelinks]{hyperref}
\hypersetup{unicode=true,pdfencoding=auto,pdftitle={A Single-Point Obstruction to Local Minimality in the P\'olya--Szeg\H{o} Capacity--Area Problem},pdfauthor={Rubin Zhang and Ning Zhang},pdfsubject={Newtonian capacity and convex geometry},pdfkeywords={Newtonian capacity, Polya-Szego conjecture, convex bodies, Alexandrov points, curvature measures}}
\usepackage[a4paper,margin=31mm,top=27mm,bottom=30mm]{geometry}
\allowdisplaybreaks[2]
\numberwithin{equation}{section}
\newtheorem{theorem}{Theorem}[section]
\newtheorem{lemma}[theorem]{Lemma}
\newtheorem{corollary}[theorem]{Corollary}
\newtheorem{proposition}[theorem]{Proposition}

\theoremstyle{remark}
\newtheorem{remark}[theorem]{Remark}
\newcommand{\R}{\mathbb R}
\newcommand{\Sph}{\mathbb S}
\newcommand{\supp}{\operatorname{supp}}
\newcommand{\Id}{\operatorname{Id}}

\newcommand{\cH}{\mathcal H}
\newcommand{\dd}{\,\mathrm d}
\newcommand{\Om}{\Omega}
\newcommand{\cD}{\mathcal D}
\title[Single-point obstruction in the capacity--area problem]{A Single-Point Obstruction to Local Minimality\\ in the P\'olya--Szeg\H{o} Capacity--Area Problem}
\author{Rubin Zhang}
\address{School of Mathematics and Statistics, Huazhong University of Science and Technology, 1037 Luoyu Road, Wuhan, Hubei 430074, China}
\curraddr{Department of Mathematical and Statistical Sciences, University of Alberta, Edmonton, Alberta T6G 2G1, Canada}
\email{rbzhang3@163.com}
\author{Ning Zhang}
\address{School of Mathematics and Statistics, Huazhong University of Science and Technology, 1037 Luoyu Road, Wuhan, Hubei 430074, China}
\email{nzhang2@hust.edu.cn}
\thanks{Corresponding author: Ning Zhang (\texttt{nzhang2@hust.edu.cn}).}
\begin{document}
\begin{abstract}
Let $Q$ be the quotient of squared Newtonian capacity by surface area
for convex bodies in $\R^3$. Suppose the lower boundary graph has a quadratic upper bound at one
point. Assume also that, along a sequence $r_j\downarrow0$, its
quadratic-scale blow-up is $H\circ\gamma$. Here $\gamma$ is a positive
convex gauge, and $H$ is a positive, nondecreasing convex profile. We redistribute the curvature measure $D^2H$ and use supporting
planes to build convex competitors. Two cancellations give a capacity change of $o(r_j^4)$ and
an area gain of $\kappa r_j^4+o(r_j^4)$ with $\kappa>0$. The competitors
converge to the original body in Hausdorff distance and have strictly
smaller $Q$ for large $j$. The gauge may be nonsmooth and $D^2H$ may be singular.
The criterion holds if there is a single Alexandrov point with a
positive-definite Hessian. For any full-dimensional local minimizer, the spatial Gaussian
curvature measure is singular with respect to surface area and the
Alexandrov Hessian has rank at most one wherever it exists. This excludes
boundaries covered by convex $W^{2,2}$ graph charts, including $C^{1,1}$
boundaries. An example with a dyadically self-similar lower graph has purely atomic
radial curvature. It gives the same descent without any positive-definite
Alexandrov point.
\end{abstract}
\maketitle

\section{Introduction}\label{sec:intro}
\subsection{The capacity--area conjecture}
For a compact set $K\subset\R^3$ of positive Newtonian capacity, we normalize
\begin{equation}\label{eq:capdef}
 C(K)=\frac1{4\pi}\inf\left\{\int_{\R^3}|\nabla v|^2\dd X:
 v\in\cD^{1,2}(\R^3),\quad v\geq1\ \text{quasi-everywhere on }K\right\},
\end{equation}
so that $C(B_R)=R$. For a full-dimensional convex body put
$S(K)=\cH^2(\partial K)$ and
\begin{equation}\label{eq:Q}
 Q(K)=\frac{C(K)^2}{S(K)}.
\end{equation}
The quotient does not change under dilations. For a planar convex plate,
we take the surface area to be twice its planar area. In this convention a planar disk
$\mathbb D_R$ of radius $R$ has
\[
 C(\mathbb D_R)=\frac{2R}{\pi},\qquad
 S(\mathbb D_R)=2\pi R^2,\qquad Q(\mathbb D_R)=\frac2{\pi^3}.
\]
The three-dimensional P\'olya--Szeg\H{o} conjecture asserts that
\begin{equation}\label{eq:PSconjecture}
 Q(K)\geq\frac2{\pi^3}
\end{equation}
for all compact convex sets of affine dimension at least two, including planar
plates, with equality only for planar disks
\cite{PolyaSzego,CrastaFragalaGazzola,PolyaSzegoBook}.
Fusco and Zhong \cite[Conjecture~1.4]{FZ} state the equivalent
conjecture that no global minimizer has nonempty interior. Here we study
local boundary conditions that
allow $Q$ to decrease among full-dimensional convex bodies arbitrarily
close in Hausdorff distance.

\subsection{Earlier results}
Bucur, Fragal\`a and Lamboley \cite[Theorem~4.5]{BFL} exclude minimizers
with a relatively open $C^3$ boundary patch of positive Gaussian
curvature. Their local convexity-preserving variations also exclude
Hausdorff-local minimizers \cite[Theorem~1.5]{FZ}.
Fusco and Zhong \cite[Theorem~1.6]{FZ} remove smoothness and strict-convexity
assumptions under two one-dimensional size conditions: the singular
boundary points and the nonregular normal directions each have zero
$\mathcal H^1$ measure. They also require the union of supporting faces
with more than one contact point to be closed. More precisely, their
closedness assumption is on $\tau_K(\Sph^2\setminus\operatorname{reg}_n K)$,
where $\operatorname{reg}_n K$ is the set of normals having a unique
contact point and $\tau_K$ takes a set of normals to the union of its
supporting faces. Their local criterion also uses an open set of
regular points with uniquely exposed supporting normals
\cite[Corollary~1.8]{FZ}. Their proof uses geometric regularity and the
Euler--Lagrange equation. Our construction uses a single quadratic-scale
tangent and allows singular or flat boundary elsewhere.

Lamboley, Novruzi and Pierre \cite{LNP} derive shape-derivative estimates
for nonsmooth domains. Their geometric application in higher dimensions deals with interior
energies \cite[Corollary~4.9]{LNP}. For exterior problems in dimensions greater than two, their estimates
do not give enough Sobolev control for that argument \cite[Remark~4.10]{LNP}. We pull the Dirichlet
energy back to the fixed exterior and pass to a half-space limit in the
localized energy term (Lemmas~\ref{lem:material} and~\ref{lem:halfspace}).

For Mahler's volume product, Reisner, Sch\"utt and Werner
\cite{ReisnerSchuttWerner} used a single generalized positive-curvature
point to exclude local minimizers. In the centrally symmetric setting,
Harrell, Henrot and Lamboley \cite{HarrellHenrotLamboley} also proved a
Gauss-curvature restriction for local minimizers. Meyer and Reisner
\cite[Theorem~1]{MeyerReisner} showed that a single point of positive
generalized Gauss curvature rules out local maximizers of the isotropy
constant. For Newton's problem
of least resistance, Lachand-Robert and Peletier
\cite{LachandRobertPeletier} showed that a minimizing convex graph cannot
be strictly convex on an open set. These results concern other
functionals. Here we must control the nonlocal exterior Dirichlet energy.
Our criterion also allows nonhomogeneous radial blow-ups and singular
radial curvature measures.
Related capacity estimates concern
thin bodies \cite{FragalaGazzolaPierre}, global $p$-capacity--area
inequalities \cite{Xiao}, and curvature pinching (X.~Jin's preprint
\emph{Curvature-Based Bounds for $p$-Capacity and Capacity--Area
Inequalities}, arXiv:2610.07564, 2026). For capacitary Brunn--Minkowski theory, see
\cite{Borell,CaffarelliJerisonLieb,ColesantiSalani}. For regularity of
planar optimal convex shapes, see \cite{LNPPlane}.

\subsection{Main results}
A \emph{positive convex gauge} is a finite convex function
$\gamma:\R^2\to[0,\infty)$ with
\begin{equation}\label{eq:gauge}
 \gamma(t\xi)=t\gamma(\xi)\ (t\geq0),\qquad
 \gamma(\xi)>0\ (\xi\ne0).
\end{equation}
Its unit sublevel set $\mathscr W=\{\gamma<1\}$ is a bounded convex
neighborhood of $0$. A \emph{positive convex radial profile} is a finite,
convex, nondecreasing function $H:[0,\infty)\to[0,\infty)$ with
$H(0)=0$ and $H(t)>0$ for $t>0$. The gauge need not be symmetric or smooth.

\begin{theorem}\label{thm:main}
Let $K\subset\R^3$ be a convex body with nonempty interior. Suppose
that, after a rigid motion, $0\in\partial K$ and the boundary of $K$
is locally represented by a convex lower graph $x_3=f(x')$, with $K$
lying above the graph.

Assume that, for some $M>0$,
\begin{equation}\label{eq:upperquadratic}
 f(0)=0,\qquad 0\le f(x')\le M|x'|^2
\end{equation}
for all $x'$ sufficiently close to the origin. Also assume there are a sequence $r_j\downarrow0$, a positive convex
gauge $\gamma$, and a positive convex radial profile $H$ such that
\begin{equation}\label{eq:blowup}
 p_j(\xi):=\frac{f(r_j\xi)}{r_j^2}\longrightarrow
 P(\xi):=H(\gamma(\xi))
 \quad\text{locally uniformly in }\R^2.
\end{equation}
Then $K$ is not a local minimizer of $Q$. In fact, there are
full-dimensional convex bodies $K_j$ and a constant $\kappa>0$ such that
\begin{equation}\label{eq:rates}
 \begin{aligned}
 d_H(K_j,K)&=O(r_j^2),\\
 C(K_j)-C(K)&=o(r_j^4),\\
 S(K_j)-S(K)&=\kappa r_j^4+o(r_j^4).
 \end{aligned}
\end{equation}
The bodies $K_j$ can be chosen to coincide with $K$ outside
$B_{C_0r_j}(0)$, for a constant $C_0$ independent of $j$.
For all sufficiently large $j$, we have $Q(K_j)<Q(K)$.
\end{theorem}

The upper bound \eqref{eq:upperquadratic} is required at all small
scales, while the convergence in \eqref{eq:blowup} is needed only along
one sequence. The gauge may be nonsmooth, and $D^2H$ may be singular.

\begin{corollary}\label{cor:alex}
Let $K$ be a full-dimensional convex body. Suppose that at some boundary
point a lower graph satisfies \eqref{eq:upperquadratic} and that
$f(r_j\xi)/r_j^2$ converges locally uniformly, along a sequence
$r_j\downarrow0$, to a finite convex function $P$ such that
\[
 P(t\xi)=t^2P(\xi)\quad(t\ge0),\qquad
 P(\xi)>0\quad(\xi\ne0).
\]
Then $K$ is not a local minimizer of $Q$. In particular, this holds if
the lower graph has a positive-definite Alexandrov expansion
\begin{equation}\label{eq:Alex}
 f(x')=\frac12 x'\cdot A x'+o(|x'|^2),\qquad A>0,
\end{equation}
with $f(0)=0$ and $Df(0)=0$.
\end{corollary}
\begin{proof}
The set $\{P<1\}$ is a bounded convex neighborhood of the origin:
convexity gives convexity of this set, and continuity and positivity
of $P$ on the unit circle give its inner and outer radius bounds.
Its Minkowski gauge is $\gamma=\sqrt{P}$, by the two-homogeneity of $P$.
Thus $\gamma$ is a positive convex gauge and $P=H\circ\gamma$ with
$H(t)=t^2$. Theorem~\ref{thm:main} applies.
For the final assertion, convexity and $Df(0)=0$ give $f\ge0$ near the
origin, and \eqref{eq:Alex} gives both \eqref{eq:upperquadratic} and
the stated convergence with $P(\xi)=\xi\cdot A\xi/2$.
\end{proof}

The \emph{spatial Gaussian curvature measure} $\mathcal G_K$ assigns to a
boundary patch the spherical area of its multivalued Gauss image
(Section~\ref{sec:prelim}).
\begin{corollary}\label{cor:curvatureintro}
Let $K$ be a full-dimensional local minimizer of $Q$. Then $\mathcal G_K$
is singular with respect to $\cH^2|_{\partial K}$. At every Alexandrov
point of any local convex graph chart, its Alexandrov Hessian has rank at
most one. So the generalized second fundamental form has rank at most
one at surface-almost-every point. In particular, a convex body whose
entire boundary is locally representable by convex $W^{2,2}$ graphs, or
whose boundary is $C^{1,1}$, cannot be a local minimizer.
\end{corollary}

For smooth convex bodies with uniform inner and outer ball bounds,
Corollary~\ref{cor:energyblowup} also shows that the squared-curvature
energy must diverge along any sequence converging to a full-dimensional
local minimizer. Corollary~\ref{cor:alex} covers nonsmooth
positive two-homogeneous tangents such as
$P(\xi)=\max\{\xi_1^2,\xi_2^2\}$. The next example uses the same square
gauge with a nonhomogeneous profile.
\begin{theorem}\label{thm:dyadicintro}
Define $H_{\mathrm d}(0)=0$ and, for each $n\in\mathbb Z$,
\begin{equation}\label{eq:Hdyadicintro}
 H_{\mathrm d}(t)=3\cdot2^nt-2\cdot4^n,
 \qquad 2^n\le t\le 2^{n+1}.
\end{equation}
Let $\gamma_{\mathrm{sq}}(x_1,x_2)=\|(x_1,x_2)\|_\infty=\max\{|x_1|,|x_2|\}$,
$f_{\mathrm d}=H_{\mathrm d}\circ\gamma_{\mathrm{sq}}$, and
\begin{equation}\label{eq:Kdyadicintro}
 K_{\mathrm d}=\bigl\{(x',x_3):|x_1|,|x_2|\le1,
 f_{\mathrm d}(x')\le x_3\le2\bigr\}.
\end{equation}
No boundary point of $K_{\mathrm d}$ has a positive-definite Alexandrov
second derivative, while Theorem~\ref{thm:main} applies at the origin.
For $r_j=2^{-j}$ there are explicit convex competitors satisfying
\begin{equation}\label{eq:explicitrateintro}
 d_H(K_j,K_{\mathrm d})\le\tfrac14r_j^2,\qquad
 C(K_j)-C(K_{\mathrm d})=o(r_j^4),\qquad
 S(K_j)-S(K_{\mathrm d})=\tfrac{63}{64}r_j^4+O(r_j^6).
\end{equation}
\end{theorem}

Section~\ref{sec:example} shows why the nonregular edges and flat
supporting contacts prevent direct use of the criteria in
\cite[Theorem~1.6 and Corollary~1.8]{FZ}.

\subsection{The two cancellations}
Redistributing the radial curvature $D^2H$ gives a
displacement $\Phi=\psi\circ\gamma$ with compact support. For small $\eta>0$,
$P+\eta\Phi$ is convex, and the displacement satisfies
\[
 \int_{\R^2}\Phi\dd\xi=0,
 \qquad
 \int_{\R^2}DP\cdot D\Phi\dd\xi=0.
\]
The first identity cancels the leading capacity term in the half-space
limit. The second cancels the mixed term in the area expansion. The remaining
positive term is proportional to $\int_{\R^2}|D\Phi|^2\dd\xi$.
For the profile $H(t)=t^2$, the two identities are equivalent. For a
general profile, we impose both conditions separately. We use supporting planes to move the perturbation to the original
graph. A patch of diameter $O(r)$ and
displacement height $O(r^2)$ then gives a positive area increment of
order $r^4$. The fixed-exterior material-energy error is $O(r^5)$ because
the local equilibrium energy is $O(r^3)$ and the deformation gradient
is $O(r)$. These cancellations leave a capacity increment of $o(r^4)$.

\subsection{Scope}
The criterion applies to subsequential tangents with homothetic gauge
level sets. A quadratic tangent with a zero-curvature direction can prevent
compact convex replacement by this method; see
Section~\ref{sec:degenerate}. The global P\'olya--Szeg\H{o} conjecture and
the existence of full-dimensional local minimizers are outside its scope.

\section{Preliminaries and notation}\label{sec:prelim}
\subsection{Convex bodies, Hausdorff convergence, and graph charts}
Let $\mathcal K^3_0$ denote the compact convex subsets of $\R^3$ with
nonempty interior. A body $K\in\mathcal K^3_0$ is a \emph{local minimizer}
of $Q$ if there exists $\varepsilon>0$ such that
$Q(K)\le Q(L)$ for every $L\in\mathcal K^3_0$ with
$d_H(K,L)<\varepsilon$. We do not fix the scale of the competitors.
For nonempty compact sets $K,L$ their Hausdorff distance is
\[
 d_H(K,L)=\max\left\{\sup_{x\in K}\operatorname{dist}(x,L),
                    \sup_{y\in L}\operatorname{dist}(y,K)\right\}.
\]
For compact convex sets, the support functions $h_K(e)=\sup_{x\in K}x\cdot e$
give $d_H(K,L)=\sup_{e\in\Sph^2}|h_K(e)-h_L(e)|$. If $K_j\to K$ in this distance and $K$ has nonempty interior, then
fixed inner and outer balls bound $K_j$ for all sufficiently large $j$. Volume and surface area are
continuous under this convergence \cite[Chapters~1--4]{Schneider}.

A \emph{lower boundary graph} of $K$ is a chart found after a rigid
motion. Near the chosen point, the boundary agrees with
$\{(x',f(x'))\}$ for a finite, locally Lipschitz convex function $f$, with
$K$ locally above its graph. Here $x'=(x_1,x_2)\in\R^2$. Denote the
horizontal projection by $\Pi_K$. Every vertical fiber is an interval.
On $\Pi_K^\circ$, the lower endpoint $\varphi_-$ is convex and the upper
endpoint $\varphi_+$ is concave. Lemma~\ref{lem:realize} uses these fibers
to realize a local graph replacement as a convex body. The area of a
Lipschitz graph over a planar set $E$ is
$\int_E\sqrt{1+|Df|^2}\dd x'$.

\subsection{Subdifferentials, Alexandrov points, and blow-up stability}
For a finite convex function $f$ on an open convex planar set $\Pi$, its
subdifferential is
\[
 \partial f(x')=\{p\in\R^2:f(y')\ge f(x')+p\cdot(y'-x')
                   \text{ for all }y'\in\Pi\}.
\]
Each $p\in\partial f(x')$ determines a supporting affine function
$y'\mapsto f(x')+p\cdot(y'-x')$. Nearby secant slopes bound subgradients
on compact subsets of $\Pi$. If finite convex functions satisfy $f_j\to f$ locally uniformly,
then the following hold.
\begin{enumerate}[label=\textup{(\roman*)},leftmargin=2.5em]
\item if $x_j\to x$ and $p_j\in\partial f_j(x_j)$ with $p_j\to p$,
then $p\in\partial f(x)$.
\item at almost every point at which $f$ is differentiable,
$Df_j\to Df$ (taking gradients where they exist).
\item $f_j\to f$ strongly in $W^{1,q}$ on smaller compact subsets for
each finite $q$, provided the functions are considered on a common
slightly larger neighborhood.
\end{enumerate}
Taking limits in the supporting inequality gives \textup{(i)}.
At each point where $f$ is differentiable and the gradients are defined,
the bounded gradients of $f_j$ can only approach $Df$. This gives
\textup{(ii)}. Dominated convergence gives \textup{(iii)}. See \cite[Chapter~12]{RockafellarWets}
for the corresponding stability results for subdifferentials.

A point $x_0$ is an \emph{Alexandrov point} of $f$ if $f$ is differentiable
there and some symmetric matrix $D_A^2 f(x_0)$ satisfies
\[
 f(x)=f(x_0)+Df(x_0)\cdot(x-x_0)
   +\tfrac12(x-x_0)\cdot D_A^2f(x_0)(x-x_0)
   +o(|x-x_0|^2).
\]
For a convex function, $D_A^2f(x_0)$ is positive semidefinite and exists
at almost every point (Alexandrov's theorem)
\cite{EvansGariepy,FigalliBook}. At a point with a unique tangent plane,
rotating that plane to horizontal changes a positive-definite quadratic
tangent by an invertible congruence and a positive scalar factor.
The \emph{generalized second fundamental form} at an Alexandrov point
is the second fundamental form of the resulting quadratic graph
expansion, with first fundamental metric $I+Df\otimes Df$.
Its rank equals $\operatorname{rank}D_A^2f$.

A positive convex gauge $\gamma$ is globally Lipschitz and satisfies
$m_\gamma|\xi|\le\gamma(\xi)\le M_\gamma|\xi|$ for some
$m_\gamma,M_\gamma>0$. Its sublevel sets
$\{\gamma<t\}$ are homothetic and may be asymmetric \cite{Schneider}.

\subsection{Monge--Amp\`ere and spatial curvature measures}
For a convex graph $f$ over $\Pi$, its Alexandrov Monge--Amp\`ere
measure is $M_f(E)=\mathcal L^2(\partial f(E))$ for Borel $E\subset\Pi$.
Under locally uniform convergence of convex functions these measures
converge weakly on compact subdomains, and the absolutely continuous
part of $M_f$ has density $\det D_A^2 f$ with respect to planar Lebesgue
measure \cite{DePhilippisFigalli,FigalliBook}. A singular part can remain even when the convex function is
continuously differentiable.

The \emph{spatial Gaussian curvature measure} $\mathcal G_K$ assigns to a
boundary patch the spherical area of its generalized Gauss image, using
all outward unit normals at nonsmooth points. Its total mass is $4\pi$
\cite{Federer,Schneider}. The slope-to-normal map for a lower graph is
$N(p)=(p,-1)/(1+|p|^2)^{1/2}$, with spherical area Jacobian
$(1+|p|^2)^{-3/2}$. On the base plane, the curvature measure is
\begin{equation}\label{eq:Gaussimage}
 \mathcal G_f(E)=\int_{\partial f(E)}\frac{\dd p}{(1+|p|^2)^{3/2}}.
\end{equation}
The weight in \eqref{eq:Gaussimage} is at most one, so
\[
 \mathcal G_f(E)\le M_f(E)\qquad\text{for every Borel }E\subset\Pi.
\]
So if $M_f$ is absolutely continuous with respect to planar Lebesgue
measure, then $\mathcal G_f$ is also absolutely continuous. The absolutely
continuous part of $\mathcal G_f$ has densities
\begin{equation}\label{eq:Gaussdensity}
 \frac{\dd\mathcal G_f^{\mathrm{ac}}}{\dd x'}
 =\frac{\det D_A^2 f}{(1+|Df|^2)^{3/2}},\qquad
 \frac{\dd\mathcal G_f^{\mathrm{ac}}}{\dd\mathcal H^2|_{\operatorname{graph}f}}
 =\frac{\det D_A^2 f}{(1+|Df|^2)^2}.
\end{equation}
Dividing the planar density by the graph area factor $\sqrt{1+|Df|^2}$
gives the surface density. The full measure may also have singular mass at corners, edges, or
other nonsmooth parts of the boundary
\cite{Federer,Schneider,FigalliBook}.

\subsection{Newtonian capacity on a fixed exterior}
The homogeneous Sobolev space $\cD^{1,2}(\R^3)$ is the completion of
$C_c^\infty(\R^3)$ in the gradient norm, represented in $L^6$.
The minimizer in \eqref{eq:capdef}, denoted by $u=u_K$, is unique,
satisfies $0\le u\le1$, equals one quasi-everywhere on $K$, is harmonic
in $\Om=\R^3\setminus K$, and tends to zero at infinity. For a
full-dimensional convex body its exterior Dirichlet boundary is regular,
so the harmonic potential has boundary value one in the usual continuous
sense \cite{Jerison}. In weak form,
\[
 \int_\Om\nabla u\cdot\nabla\varphi\dd X=0
 \qquad(\varphi\in\cD^{1,2}_0(\Om)).
\]
Here $\cD^{1,2}_0(\Om)$ denotes the zero Sobolev-trace space, whose
elements extend by zero through $K$ into $\cD^{1,2}(\R^3)$.
For the boundary blow-up we write $U=1-u$ in $\Om$ and extend $U$ by zero
in $K$. Capacity is monotone under inclusion and satisfies $C(tK)=tC(K)$
for $t>0$. Both $C(K)^2$ and $S(K)$ scale quadratically.

For a bi-Lipschitz map $F=\Id+V$ with $V$ compactly supported, pull back
the Dirichlet problem on $\R^3\setminus F(K)$ to $\Om$. The coefficient
field depends on $DF$, and the correction to the potential has zero trace
on $\partial\Om$ \cite{ElcratMiller,HenrotPierre,LNP}.
Lemma~\ref{lem:material} gives a bound for the energy change in this
formulation.

\subsection{Scaling and notation conventions}
We use $D$ for the planar gradient, $\nabla$ for the spatial gradient,
and $D^2H$ for the one-dimensional distributional second derivative.
We use $x'$ for the original planar coordinates, $\xi$ for the
rescaled planar coordinates, and $X$ for a spatial variable. In
Section~\ref{sec:moments}, $\nu$ denotes a signed radial curvature measure.
We write $n_K$ for the outward unit normal. The notation $C(K)$
always means capacity, while a bare $C$ denotes a positive constant.
We use $h_K$ for the support function, $h=H+\eta\psi$ for the modified
radial profile, and $h_r(\xi)=f(r\xi)/r$ for the rescaled graph height.
The constants $c,C,R_*$ may change from one estimate to another.

The bound $0\le f(x')\le M|x'|^2$ gives $Df(0)=0$. For $|x'|\le Rr$, we bound each supporting slope using the secant
slopes at distance $r$ on both sides of $x'$. The quadratic bound on the larger ball of
radius $(R+1)r$ gives $\sup_{|x'|\le Rr}|Df(x')|=O_R(r)$ in the
subgradient sense. So $f(r\xi)/r$ has height and Lipschitz constant
$O_R(r)$ on $|\xi|\le R$, as used in Lemma~\ref{lem:halfspace}.

\section{Model perturbations and convex gluing}\label{sec:moments}
\subsection{Curvature redistribution and gauge-layer cancellation}
Passing to the limit in $0\le p_j\le M|\xi|^2$ gives
\[
 0\le P(\xi)\le M|\xi|^2\qquad(\xi\in\R^2).
\]
Take $e$ with $\gamma(e)=1$. The bound $0\le H(t)=P(te)\le M|e|^2t^2$ gives
\[
 H'_+(0)=0.
\]
The distributional second derivative $\mu=D^2H$ is a nonnegative
Radon measure on $(0,\infty)$.

\begin{lemma}\label{lem:blocks}
For every $\delta>0$, $\mu((0,\delta))>0$. So one may choose seven pairwise separated compact intervals $E_1,\ldots,E_7\subset(0,\infty)$, ordered from left to right, with $0<\mu(E_i)<\infty$.
\end{lemma}
\begin{proof}
If $\mu((0,\delta))=0$, then $H$ is affine on $(0,\delta)$.
Because $H(0)=H'_+(0)=0$, it must be zero there. This contradicts
$H(t)>0$ for $t>0$. So the support of $\mu$ has positive points as
close to zero as we wish. Choose small, disjoint compact intervals
around seven such points.
\end{proof}

Write $a_i=\min E_i$ and $b_i=\max E_i$, so $b_i<a_{i+1}$ for $i<7$.
Fix $0<a_0<a_1$ and put
\[
 \mathcal F_H(s)=\int_{a_0}^s tH'(t)\dd t.
\]
Because $H$ is convex, its derivative $H'$ exists almost everywhere
and is locally bounded on $(0,\infty)$.

\begin{lemma}\label{lem:moments}
There is a nonzero finite signed measure $\nu=\alpha\mu$ supported on $E_1\cup\cdots\cup E_5$, with $|\alpha|\le1$, such that
\begin{equation}\label{eq:fourmoments}
 \int\dd\nu=\int s\dd\nu(s)=\int s^3\dd\nu(s)
 =\int\mathcal F_H(s)\dd\nu(s)=0.
\end{equation}
If
\[
 \psi(t)=\int(t-s)_+\dd\nu(s),
\]
then $\psi\not\equiv0$ is compactly supported and Lipschitz, and
\begin{equation}\label{eq:twomoments}
 \int_0^\infty t\psi(t)\dd t=0,
 \qquad\int_0^\infty tH'(t)\psi'(t)\dd t=0.
\end{equation}
For $0<\eta\le\frac12$, $h=H+\eta\psi$ is convex and nondecreasing.
\end{lemma}
\begin{proof}
Set $\nu=\sum_{i=1}^5c_i\,\mu|_{E_i}$. The four homogeneous equations in \eqref{eq:fourmoments} have a nonzero solution in the five coefficients $c_i$. Rescale the coefficients so that $\max_i|c_i|\le1$. The intervals
do not overlap, and each has positive $\mu$-measure. So $\nu\ne0$.

The integrand $(t-s)_+$ is $1$-Lipschitz in $t$, so $\psi$ is Lipschitz. It vanishes for $t\le a_1$. It is also zero for $t\ge b_5$ because of the first two conditions
in \eqref{eq:fourmoments}. Since $D^2\psi=\nu\ne0$, the function
$\psi$ is not identically zero. Take $B>b_5$. Fubini's theorem gives
\begin{align*}
 \int_0^\infty t\psi(t)\dd t
 &=\int\left(\frac{B^3}{3}-\frac{B^2s}{2}+\frac{s^3}{6}\right)\dd\nu(s)
 =\frac16\int s^3\dd\nu(s)=0,\\
 \int_0^\infty tH'(t)\psi'(t)\dd t
 &=\int\left(\int_s^B tH'(t)\dd t\right)\dd\nu(s)
 =-\int\mathcal F_H(s)\dd\nu(s)=0.
\end{align*}
The curvature measure of $h$ satisfies
\[
 D^2h=(1+\eta\alpha)\mu\ge\tfrac12\mu.
\]
So $h$ is convex. Because it agrees with $H$ near zero, we have
$h'_+(0)=0$. Convexity then makes $h$ nondecreasing.
\end{proof}

Set
\[
 \Phi(\xi)=\psi(\gamma(\xi)),\qquad
 G(\xi)=h(\gamma(\xi))=P(\xi)+\eta\Phi(\xi).
\]
The function $G$ is convex because $h$ is convex and nondecreasing
and $\gamma$ is convex. The function $\Phi$ is Lipschitz and has
support in $\{\gamma\le b_5\}$.

\begin{lemma}\label{lem:layer}
Let $\mathscr W=\{\gamma<1\}$ and $J=\int_{\mathscr W}|D\gamma|^2\dd\xi$. Then $0<J<\infty$,
\begin{equation}\label{eq:zeromom}
 \int_{\R^2}\Phi\dd\xi=0,\qquad
 \int_{\R^2}DP\cdot D\Phi\dd\xi=0,
\end{equation}
and
\begin{equation}\label{eq:positiveenergy}
 \frac12\int_{\R^2}\bigl(|DG|^2-|DP|^2\bigr)\dd\xi
 =\eta^2J\int_0^\infty t|\psi'(t)|^2\dd t=:\kappa>0.
\end{equation}
All integrals in \eqref{eq:zeromom}--\eqref{eq:positiveenergy} are over the compact region where $G-P$ may be nonzero.
\end{lemma}
\begin{proof}
Homogeneity gives $D\gamma(t\xi)=D\gamma(\xi)$ at differentiability points for $t>0$. Scaling the sublevel sets then gives
\[
 |\{\gamma<t\}|=t^2|\mathscr W|,
 \qquad\int_{\{\gamma<t\}}|D\gamma|^2\dd\xi=t^2J.
\]
For a compactly supported function $a$ with well-defined integrals,
the pushforward measures under $\gamma$ give
\[
 \int_{\R^2}a(\gamma)\dd\xi=2|\mathscr W|\int_0^\infty ta(t)\dd t,
 \quad\int_{\R^2}a(\gamma)|D\gamma|^2\dd\xi
 =2J\int_0^\infty ta(t)\dd t.
\]
Almost everywhere, $DP=H'(\gamma)D\gamma$ and $D\Phi=\psi'(\gamma)D\gamma$. Use the first integration formula with $a=\psi$ and the second
with $a=H'\psi'$. The conditions in \eqref{eq:twomoments} give
\eqref{eq:zeromom}. Now expand $|DP+\eta D\Phi|^2$ and use the
second formula with $a=|\psi'|^2$. This gives
\eqref{eq:positiveenergy}. The radial integral is positive because $\psi$ is nonzero, Lipschitz,
and has compact support. Also, $J<\infty$ because $\gamma$ is
Lipschitz and $\mathscr W$ is bounded. If $J=0$, then $\gamma$ is constant on the open convex set $\mathscr W$. It vanishes at zero, so homogeneity forces $\gamma\equiv0$, contrary to \eqref{eq:gauge}.
\end{proof}

\begin{remark}\label{rem:two}
When $H(t)=t^2$, integration by parts gives $\int tH'(t)\psi'(t)\dd t=-4\int t\psi(t)\dd t$, so the two cancellation conditions in \eqref{eq:twomoments} are equivalent. For a nonhomogeneous profile, the two conditions are generally independent.
Lemma~\ref{lem:moments} gives both conditions for the measure $D^2H$, even
when this measure is singular.
\end{remark}

\subsection{Supporting-plane gaps and stable gluing}\label{sec:gluing}
For $y\ne0$, put $s=\gamma(y)$ and let $q\in\partial P(y)$. Apply the supporting inequality along the ray through $y$ to get
$\beta=q\cdot y/s\in\partial H(s)$, with $\beta\ge0$. For $x\ne0$, the point $sx/\gamma(x)$ has the same $P$-value as $y$. The supporting inequality at this point gives $q\cdot x\le\beta\gamma(x)$, which also holds at $x=0$. With $t=\gamma(x)$, we obtain
\[
 P(x)-P(y)-q\cdot(x-y)
 \ge H(t)-H(s)-\beta(t-s).
\]
For a one-dimensional convex function, integration of its second derivative measure gives
\begin{align*}
 H(t)-H(s)-\beta(t-s)&\ge\int_{(t,s)}(v-t)\dd\mu(v),&&t<s,\\
 H(t)-H(s)-\beta(t-s)&\ge\int_{(s,t)}(t-v)\dd\mu(v),&&t>s.
\end{align*}
These bounds hold for every $\beta\in\partial H(s)$, even if
$\mu$ has atoms at the endpoints. The same inequalities hold for $G=h\circ\gamma$ with $D^2h$ in place of $\mu$.

At $y=0$, the quadratic bound on $P$ and the equality $G=P$ near zero
show that both gradients are zero. So the supporting plane of each
function at the origin is zero. Since $h(0)=h'_+(0)=0$, for $t=\gamma(x)>0$ we have
\[
 G(x)=h(t)=\int_{(0,t)}(t-v)\,\dd D^2h(v).
\]
The perturbation in Lemma~\ref{lem:moments} leaves the curvature on $E_6$ and $E_7$ unchanged. Choose
\[
 b_5<R_1<a_6\le b_6<R_2<a_7\le b_7<R_3<R_4,
\]
and fix $\varepsilon>0$ with $b_7<R_3-\varepsilon<R_3+\varepsilon<R_4$.
For $\gamma(y)\in[R_2,R_4]$ and $\gamma(x)\le R_1$, the contribution of $E_6$ bounds the gap between $P(x)$ and any supporting plane of $P$ at $y$ from below by
\begin{equation}\label{eq:deltain}
 \delta_{\rm in}:=(a_6-R_1)\mu(E_6)>0.
\end{equation}
For $\gamma(y)\le R_2$ and $R_3-\varepsilon\le\gamma(x)\le R_3+\varepsilon$, the contribution of $E_7$ bounds the gap between $G(x)=P(x)$ and any supporting plane of $G$ at $y$ from below by
\begin{equation}\label{eq:deltaout}
 \delta_{\rm out}:=(R_3-\varepsilon-b_7)\mu(E_7)>0.
\end{equation}
At $y=0$, the formula for $G$ gives
\[
 G(x)-0\ge\int_{E_7}(\gamma(x)-v)\,\dd\mu(v)
 \ge\delta_{\rm out}.
\]
Choose $\eta>0$ in Lemma~\ref{lem:moments} small enough that
$\eta\|\psi\|_\infty<\delta_{\rm in}/4$.

\begin{lemma}\label{lem:gluing}
Let $p_j$ be convex functions on expanding neighborhoods of the origin, with $p_j\to P$ locally uniformly, where $P,G$ are as constructed above. There exist convex functions $G_j$ on the same large neighborhoods such that
\begin{enumerate}[label=\textup{(\roman*)},leftmargin=2.5em]
\item $G_j=p_j$ outside one fixed compact set $\{\gamma\le R_3\}$;
\item $G_j=G$ on $\{\gamma\le R_1\}$ for all sufficiently large $j$;
\item $G_j-p_j\to G-P=\eta\Phi$ uniformly and strongly in $H^1$ on every fixed region containing the support, and $\|D(G_j-p_j)\|_{L^\infty}$ is uniformly bounded.
\end{enumerate}
\end{lemma}
\begin{proof}
Let $T_j$ be the supremum of the supporting planes of $G$ at points in $\{\gamma\le R_2\}$ and of $p_j$ at points in $\{R_2\le\gamma\le R_4\}$. The contact points lie in fixed compact sets. Local secant bounds
give a common bound on the subgradients. So $T_j$ is finite and convex.

On $\{\gamma\le R_1\}$, the planes of $G$ give $T_j\ge G$. For large $j$, every plane of $p_j$ in the supremum lies at least $\delta_{\rm in}/2$ below $P$ there. Otherwise, along a subsequence there are points $x_j$ with $\gamma(x_j)\le R_1$, contact points $y_j$ with $R_2\le\gamma(y_j)\le R_4$, and $q_j\in\partial p_j(y_j)$ such that
\[
 p_j(y_j)+q_j\cdot(x_j-y_j)>P(x_j)-\delta_{\rm in}/2.
\]
Pass to a further subsequence so that $x_j\to x$, $y_j\to y$,
and $q_j\to q$. For each fixed $z$, take the limit in the supporting
inequality $p_j(z)\ge p_j(y_j)+q_j\cdot(z-y_j)$. This gives
$q\in\partial P(y)$, as in Section~\ref{sec:prelim}. The limiting plane lies at most $\delta_{\rm in}/2$ below $P(x)$.
This contradicts \eqref{eq:deltain}. Since $\|G-P\|_\infty<\delta_{\rm in}/4$, the planes of $p_j$ lie below $G$ there, and $T_j=G$.

On the collar $R_3-\varepsilon\le\gamma\le R_3+\varepsilon$, the planes of $p_j$ give $T_j\ge p_j$. By \eqref{eq:deltaout} and local uniform convergence, the planes of $G$ lie below $p_j$ there for large $j$. So $T_j=p_j$ throughout the collar.

On $R_1\le\gamma\le R_3$ we have $G=P$, since $b_5<R_1$. The planes in the supremum give $T_j\le\max\{P,p_j\}$. They also give $T_j\ge P$ for $\gamma\le R_2$ and $T_j\ge p_j$ for $\gamma\ge R_2$ in this annulus. So $T_j\to G$ uniformly on $\{\gamma\le R_3\}$. Set $G_j=T_j$ on $\{\gamma<R_3\}$ and $G_j=p_j$ outside. The two functions agree near $\{\gamma=R_3\}$. So $G_j$ is
locally convex and is convex on the original convex domain.

The patched functions converge locally uniformly to $G$ because
$G=P$ outside $\{\gamma<R_3\}$. Local uniform convergence of convex functions gives almost-everywhere gradient convergence. Secant bounds on a slightly larger fixed set bound the gradients of $G_j$ and $p_j$ uniformly. Dominated convergence gives strong $H^1$ convergence of $G_j$ to $G$ and of $p_j$ to $P$ on the fixed set containing the support of the displacement. Subtracting and using (i) gives (iii).
\end{proof}

\section{Capacitary estimates}\label{sec:potential}
\subsection{A fixed-exterior material estimate}
Let $u$ be the equilibrium potential of $K$ and write $\Om=\R^3\setminus K$. For a compactly supported Lipschitz vector field $V$, set
\[
 B_V=(\operatorname{div}V)I-DV-DV^T.
\]
\begin{lemma}\label{lem:material}
There are absolute constants $\eta_0,c>0$ such that, whenever $V$ is
compactly supported and Lipschitz with $\|DV\|_\infty\le\eta_0$,
$F=\Id+V$ is bi-Lipschitz and
\begin{equation}\label{eq:material}
 \left|4\pi\bigl(C(F(K))-C(K)\bigr)-
 \int_\Om B_V\nabla u\cdot\nabla u\dd X\right|
 \le c\int_\Om |DV|^2|\nabla u|^2\dd X.
\end{equation}
\end{lemma}
\begin{proof}
If $\|DV\|_\infty<1$, the Lipschitz bound for $V$ gives a
two-sided Lipschitz bound for $F$. Also, $X=Y-V(X)$ has a unique
solution for every $Y$ by the contraction principle. The map $F$ equals the identity outside a compact set. Pulling the energy on $F(\Om)$ back to $\Om$ gives the coefficient matrix
\[
 A_V=\det(I+DV)(I+DV)^{-1}(I+DV)^{-T}.
\]
Expanding at $DV=0$ gives
\[
 A_V=I+B_V+R_V,\qquad |R_V|\le c|DV|^2.
\]
Write the pulled-back potential as $u+w$, where $w$ has zero Sobolev trace on $\partial\Om$. Subtracting the weak equation for $u$ from the equation for $u+w$ gives
\[
 \int_\Om A_V\nabla w\cdot\nabla\varphi
 =-\int_\Om(A_V-I)\nabla u\cdot\nabla\varphi
 \quad(\varphi\in\cD^{1,2}_0(\Om)).
\]
Test with $w$ and use uniform ellipticity to obtain
\[
 \|\nabla w\|_{L^2(\Om)}\le
 c\bigl\|\,|DV|\,|\nabla u|\,\bigr\|_{L^2(\Om)}.
\]
Test the equation for $u+w$ with $w$. The cross term in the energy
then vanishes, and we get
\[
 4\pi C(F(K))=
 \int_\Om A_V\nabla u\cdot\nabla u
 -\int_\Om A_V\nabla w\cdot\nabla w.
\]
Subtract $4\pi C(K)=\int_\Om |\nabla u|^2$. The bound on $R_V$ controls the coefficient remainder, and the estimate for $\nabla w$ controls the last integral. Both are bounded by a constant times
$\int_\Om |DV|^2|\nabla u|^2$. This proves \eqref{eq:material}.
\end{proof}

\subsection{Half-space blow-up of the energy}
\begin{lemma}\label{lem:halfspace}
Assume $0\le f(x')\le M|x'|^2$ near the origin. For $r>0$, put
$h_r(\xi)=f(r\xi)/r$ and, for $t<0$, set
\[
 w_r(\xi,t)=r^{-1}U\bigl(r\xi,r(t+h_r(\xi))\bigr),
\]
where $U=1-u$ in the exterior and zero in $K$.
For every sequence $r_j\downarrow0$ there is a subsequence, not
relabeled, and $c_*\ge0$ such that, for every fixed $R>0$,
\begin{equation}\label{eq:halfspace}
 w_{r_j}\longrightarrow -c_*t
 \quad\text{strongly in }H^1(B_R^-),
 \qquad B_R^-:=B_R\cap\{t<0\}.
\end{equation}
The traces of $w_{r_j}$ on the flat boundary $\{t=0\}$ are zero.
We also have, for each fixed $R>0$,
\begin{equation}\label{eq:energyball}
 \int_{B_{Rr}\setminus K}|\nabla u|^2\dd X=O_R(r^3)
 \quad\text{as }r\downarrow0.
\end{equation}
\end{lemma}
\begin{proof}
The secant bound in Section~\ref{sec:prelim} gives $\sup_{|x'|\le Rr}|Df(x')|=O_R(r)$ in the subgradient sense. For small $a$, the ball $B_a(ae_3)$ lies in the lower graph chart and is internally tangent to $K$ at zero. To see this, its lower hemisphere satisfies $a-\sqrt{a^2-|x'|^2}\ge|x'|^2/(2a)\ge f(x')$, and the ball stays within the neighborhood in which $K$ lies above the graph.

Comparison with the potential of this ball gives, for $X\in\Om$,
\[
 0\le U(X)\le1-\frac{a}{|X-ae_3|}.
\]
The rescaled graph $h_r(X')=f(rX')/r$ has height and slope $O_R(r)$ on fixed bounded sets. Flatten it by $(X',t)\mapsto(X',t+h_r(X'))$. On a fixed lower half-ball write
\[
 w_r(X',t)=r^{-1}U\bigl(rX',r(t+h_r(X'))\bigr),\qquad t<0.
\]
The coefficient matrices $A_r$ in the flattened equation satisfy $\|A_r-I\|_{L^\infty(B_{2R}^-)}=O_R(r)$. Choose $\theta\in C_c^\infty(B_{2R})$ equal to one on $B_R$. Since $w_r$ has zero trace on $t=0$, $\theta^2w_r$ is an admissible test function. Testing the equation gives
\[
 \int_{B_R^-}|\nabla w_r|^2\dd X
 \le C_R\int_{B_{2R}^-}|w_r|^2\dd X\le C_R.
\]
The last bound follows from the ball comparison. It bounds $w_r$ on
each fixed half-ball because $|r(X',t+h_r)-ae_3|=a+O_R(r)$. The flattening maps and their inverses have uniformly bounded derivatives. Undo the flattening and use $\nabla_X(U(rX)/r)=\nabla U(rX)$. A fixed half-ball of slightly larger radius covers the preimage of
the rescaled exterior patch. The volume factor $r^3$ then gives
\eqref{eq:energyball}.

Take a diagonal subsequence such that on each bounded half-ball $w_r\rightharpoonup w_\infty$ in $H^1$ and $w_r\to w_\infty$ in $L^2$. Because $A_r\to I$ locally uniformly, the limit is harmonic in
the lower half-space. It has zero trace on the flat boundary by weak $H^1$ convergence. Testing the difference of the weak equations with $\theta^2(w_r-w_\infty)$ gives
\[
 \int_{B_R^-}|\nabla(w_r-w_\infty)|^2\dd X
 \le C_R\Big(\|w_r-w_\infty\|_{L^2(B_{2R}^-)}^2
 +\|A_r-I\|_{L^\infty(B_{2R}^-)}^2
  \|\nabla w_r\|_{L^2(B_{2R}^-)}^2\Big)\longrightarrow0.
\]
This proves strong $H^1$ convergence up to the flat boundary. In the ball comparison, $|r(X',t+h_r)-ae_3|=a-rt+O_R(r^2)$, so $0\le w_\infty(\xi,t)\le -t/a$. Reflect it oddly across $t=0$. The reflected function is harmonic
on the whole space and grows at most linearly. By Liouville's theorem,
it is affine. The zero trace and nonnegativity for $t<0$ then give
$w_\infty(\xi,t)=-c_*t$ with $c_*\ge0$.
\end{proof}

\section{Proof of Theorem~\ref{thm:main}}\label{sec:mainproof}
\subsection{Realization as a convex body}
\begin{lemma}\label{lem:realize}
Assume the geometric conditions of Theorem~\ref{thm:main}, and set
\[
 F_j(x')=r_j^2G_j(x'/r_j).
\]
For all sufficiently large $j$, there is a convex body $K_j$ whose
lower graph is $F_j$ on the modified patch. It agrees with $K$ outside
a neighborhood of the origin of diameter $O(r_j)$ and satisfies
$d_H(K_j,K)=O(r_j^2)$.
\end{lemma}
\begin{proof}
On the boundary chart, $F_j$ is convex and equals $f$ outside a disk of radius $O(r_j)$. Lemma~\ref{lem:gluing} gives the bounds
\begin{equation}\label{eq:replacementbounds}
 \|F_j-f\|_\infty\le Cr_j^2,
 \qquad\|D(F_j-f)\|_\infty\le Cr_j.
\end{equation}
Let $\Pi_K$ be the horizontal projection of $K$. The lower graph is
defined in a neighborhood of $0$, so $0\in\Pi_K^\circ$. Write
\[
 \varphi_-(x')=\min\{x_3:(x',x_3)\in K\},\qquad
 \varphi_+(x')=\max\{x_3:(x',x_3)\in K\},\quad x'\in\Pi_K.
\]
The lower profile $\varphi_-$ is convex and agrees with $f$ near $0$.
The upper profile $\varphi_+$ is concave. Define $\varphi_{-,j}$
by gluing $F_j$ to $\varphi_-$ in the open collar where they agree.
This function is locally convex in $\Pi_K^\circ$. It is therefore
convex on every segment in $\Pi_K^\circ$.

The patch is compactly contained in $\Pi_K^\circ$, so
$\varphi_{-,j}=\varphi_-$ in a relative neighborhood of
$\partial\Pi_K$. Fix $a\in\Pi_K^\circ$. For $z\in\partial\Pi_K$,
the points $z_\tau=(1-\tau)z+\tau a$ lie in this neighborhood
for small $\tau>0$. Lower semicontinuity and convexity of $\varphi_-$ give
\[
 \varphi_-(z)\le\liminf_{\tau\downarrow0}\varphi_-(z_\tau)
 \le\limsup_{\tau\downarrow0}\varphi_-(z_\tau)
 \le\varphi_-(z).
\]
Set $\varphi_{-,j}(z)=\varphi_-(z)$ on $\partial\Pi_K$. The functions agree near the boundary, so $\varphi_{-,j}$ is lower
semicontinuous there. The limits along $z_\tau$ also let us extend
the convexity inequality to the boundary. To see this, apply it to $(1-\tau)x+\tau a$ and
$(1-\tau)y+\tau a$ for $x,y\in\Pi_K$ and let $\tau\downarrow0$.

Near $0$ the fiber height $\varphi_+-\varphi_-$ has a positive lower
bound. Since $\|\varphi_{-,j}-\varphi_-\|_\infty=O(r_j^2)$ by
\eqref{eq:replacementbounds}, $\varphi_{-,j}\le\varphi_+$ for large $j$.
The set
\[
 K_j=\{(x',x_3):x'\in\Pi_K,\quad
             \varphi_{-,j}(x')\le x_3\le\varphi_+(x')\}
\]
is compact, convex and full-dimensional. It differs from $K$ only
between the two lower graphs on the patch. Each lower endpoint changes by at most $O(r_j^2)$. So the Hausdorff
distance has the same bound. Using $(1-t)\varphi_-+t\varphi_{-,j}$ as the lower profile,
with $0\le t\le1$, gives a continuous path of convex bodies from $K$
to $K_j$.
\end{proof}

\subsection{Scale bounds and the area increment}
Take $G$, $\Phi$ and $G_j$ from Lemmas~\ref{lem:moments}--\ref{lem:gluing}
and write
\[
 \delta_j(x')=F_j(x')-f(x'),\qquad
 F_j(x')=r_j^2G_j(x'/r_j).
\]
The perturbation has support of diameter $O(r_j)$. The bounds on
the modified patch are
\begin{equation}\label{eq:mainbounds}
 \|\delta_j\|_\infty=O(r_j^2),\quad
 \|D\delta_j\|_\infty=O(r_j),\quad
 \|Df\|_\infty+\|DF_j\|_\infty=O(r_j).
\end{equation}
The rescaled difference satisfies
\[
 r_j^{-2}\delta_j(r_j\xi)=G_j(\xi)-p_j(\xi)
 \longrightarrow\eta\Phi(\xi)
\]
uniformly and strongly in $H^1$ on a fixed compact set.
Let $K_j$ be the body given by Lemma~\ref{lem:realize}.

Only the lower graph changes, so the area increment is
\[
 S(K_j)-S(K)
 =\int\left(\sqrt{1+|DF_j|^2}-\sqrt{1+|Df|^2}\right)\dd x'.
\]
The gradients are $O(r_j)$ on a planar set of area $O(r_j^2)$.
Taylor expansion and the substitution $x'=r_j\xi$ give
\begin{align}
 S(K_j)-S(K)
 &=\tfrac12\int\bigl(|DF_j|^2-|Df|^2\bigr)\dd x'+O(r_j^6)\notag\\
 &=\frac{r_j^4}{2}\int\bigl(|DG_j|^2-|Dp_j|^2\bigr)\dd\xi+O(r_j^6)\notag\\
 &=\kappa r_j^4+o(r_j^4),\label{eq:areagrowth}
\end{align}
Strong $H^1$ convergence on the rescaled patch and the positive
energy increment $\kappa$ in \eqref{eq:positiveenergy} give the
last line.

\subsection{Capacity increment and strict descent}
Choose $\chi\in C_c^\infty(\R)$ equal to one near zero, with
support in a sufficiently short fixed interval. In the boundary chart set
\[
 V_j(x',x_3)=\delta_j(x')\chi\!\left(\frac{x_3-f(x')}{r_j}\right)e_3,
\]
and extend $V_j$ by zero. The map acts on each vertical fiber
separately. For large $j$ it is strictly increasing on the fiber,
sends the lower endpoint to the new one, and fixes the upper
endpoint. Thus
\[
 (\Id+V_j)(K)=K_j.
\]
By \eqref{eq:mainbounds},
\[
 \supp V_j\subset B_{R_*r_j},\qquad
 \|V_j\|_\infty=O(r_j^2),\qquad
 \|DV_j\|_\infty=O(r_j)
\]
for a fixed $R_*>0$. Lemmas~\ref{lem:material} and \ref{lem:halfspace} give
\begin{equation}\label{eq:materror}
 \int_\Om|DV_j|^2|\nabla u|^2\dd X
 =O(r_j^2)O(r_j^3)=O(r_j^5).
\end{equation}
To find the linear term, write $h_j=h_{r_j}$ and set
\[
 \phi_j(\xi)=G_j(\xi)-p_j(\xi),\qquad
 \widetilde V_j(Y)=r_j^{-2}V_j(r_jY).
\]
For $Y=(\xi,s)$ in the rescaled chart,
\[
 \widetilde V_j(\xi,s)
  =\phi_j(\xi)\chi\bigl(s-h_j(\xi)\bigr)e_3.
\]
The fields and their derivatives are uniformly bounded on a common
compact support. Lemma~\ref{lem:gluing} gives $\phi_j\to\eta\Phi$
uniformly and in $H^1$. The gradients of the convex functions $G_j$ and $p_j$ converge
almost everywhere. So $D\phi_j\to\eta D\Phi$ almost everywhere
along the full sequence. On the support,
$\|Dh_j\|_\infty=O(r_j)$.

The map $T_j(\xi,t)=(\xi,t+h_j(\xi))$ sends the lower half-space
$t<0$ to the rescaled exterior within the chart. The matrices
$B_{\widetilde V_j}\circ T_j$ are uniformly bounded and converge
almost everywhere on the common compact support to $B_{V_\infty}$, where
\[
 V_\infty(\xi,t)=\eta\Phi(\xi)\chi(t)e_3.
\]
Set $U_j(Y)=r_j^{-1}U(r_jY)$. By the chain rule and $Dh_j\to0$,
the strong $H^1$ convergence in Lemma~\ref{lem:halfspace} gives,
on each fixed lower half-ball,
\[
 \nabla u(r_jT_j(\xi,t))=-\nabla U_j(T_j(\xi,t))
 \longrightarrow c_*e_3\quad\text{in }L^2.
\]
Using $DV_j(r_jY)=r_jD\widetilde V_j(Y)$,
$B_{V_j}(r_jY)=r_jB_{\widetilde V_j}(Y)$ and $\det DT_j=1$, along
each subsequence with the half-space limit we get
\begin{align}
 \frac1{r_j^4}\int_\Om B_{V_j}\nabla u\cdot\nabla u\dd X
 &=\int_{\{t<0\}}
 (B_{\widetilde V_j}\circ T_j)
  \nabla u(r_jT_j)\cdot\nabla u(r_jT_j)\dd\xi\dd t\notag\\
 &\longrightarrow c_*^2\int_{\R^2}\int_{-\infty}^0
    e_3\cdot B_{V_\infty}e_3\dd t\dd\xi\notag\\
 &=-c_*^2\eta\int_{\R^2}\Phi(\xi)\dd\xi=0.
 \label{eq:bulkvanish}
\end{align}
Strong $L^2$ convergence gives $L^1$ convergence of the gradient
products. The matrices are uniformly bounded, and dominated
convergence applies to the limiting product. This gives the integral
limit. To find its value,
use $e_3\cdot B_{V_\infty}e_3=-\eta\Phi(\xi)\chi'(t)$,
$\int_{-\infty}^0\chi'(t)\dd t=\chi(0)=1$, and \eqref{eq:zeromom}.

Every subsequence has a further subsequence with a half-space limit.
Equation~\eqref{eq:bulkvanish} holds for each possible $c_*$. So the
full sequence satisfies
\[
 \int_\Om B_{V_j}\nabla u\cdot\nabla u\dd X=o(r_j^4).
\]
We now use Lemma~\ref{lem:material} and \eqref{eq:materror} to get
\begin{equation}\label{eq:capacitylittleo}
 C(K_j)-C(K)=o(r_j^4).
\end{equation}

With $c_K=C(K)>0$ and $s_K=S(K)>0$, equations~\eqref{eq:areagrowth}
and \eqref{eq:capacitylittleo} give
\[
 Q(K_j)-Q(K)
 =-\frac{c_K^2}{s_K^2}\,\kappa r_j^4+o(r_j^4)<0
\]
for all sufficiently large $j$. Lemma~\ref{lem:realize} gives the localization and the Hausdorff
bound stated in Theorem~\ref{thm:main}.

\section{Gaussian curvature and Sobolev regularity consequences}\label{sec:curvature}

\begin{proof}[Proof of Corollary~\ref{cor:curvatureintro}]
At an Alexandrov point of a convex graph, let $p=Df(0)$ and
$A=D_A^2f(0)$ after translating the point to the origin. Suppose that
$A$ is positive definite. Choose orthonormal coordinates with the new
horizontal plane equal to the tangent plane and the new vertical axis
pointing into $K$. The new tangent coordinates satisfy
$y=Bx'+O(|x'|^2)$ for an invertible matrix $B$. The new height is
$g(y)=\alpha\bigl(f(x')-p\cdot x'\bigr)$, where
$\alpha=(1+|p|^2)^{-1/2}>0$. As $x'=B^{-1}y+o(|y|)$, we obtain
\[
 g(y)=\tfrac12 y\cdot\bigl(\alpha B^{-T}AB^{-1}\bigr)y+o(|y|^2).
\]
Thus the new Alexandrov Hessian is positive definite.
Corollary~\ref{cor:alex} gives a strict descent.
So a local minimizer must have $\operatorname{rank}D_A^2f\le1$
at every Alexandrov point. By Alexandrov's theorem, the determinant
is zero almost everywhere. Equation~\eqref{eq:Gaussdensity} gives
singularity of the spatial Gaussian curvature measure relative to
surface area. The graph coordinates give the same rank bound for the generalized
second fundamental form.

Suppose the boundary has a cover by convex $W^{2,2}$ graphs. On any
smaller subchart the smooth convex mollifications $f_\varepsilon$
converge strongly in $W^{2,2}$. For symmetric $2\times2$ matrices,
\[
 |\det A-\det B|\le c(|A|+|B|)|A-B|.
\]
Thus $\det D^2f_\varepsilon\to\det D^2f$ in $L^1$. Weak continuity of Monge--Amp\`ere measures gives
$M_f=(\det D^2f)\,\dd x'$.
Since $\mathcal G_f\le M_f$, the curvature measure is absolutely
continuous with respect to planar area and hence to graph surface
area. A finite set of charts gives absolute continuity on $\partial K$.
But the measure is also singular, so it must be zero. This contradicts
its total mass $4\pi$. The same conclusion holds for $C^{1,1}$ convex
graphs, which belong locally to $W^{2,2}$.
\end{proof}

\begin{remark}\label{rem:FZ}
Suppose $K$ is a local minimizer and satisfies the size assumptions
of Fusco and Zhong \cite[(3.1)]{FZ}:
\[
 \cH^1(\partial K\setminus\operatorname{reg}K)=0,\qquad
 \cH^1(\Sph^2\setminus\operatorname{reg}_n K)=0.
\]
Here $\operatorname{reg}K$ is the set of boundary points with a unique
normal, and $\operatorname{reg}_nK$ is the set of normals with a unique
contact point. Under these assumptions,
\cite[Propositions~3.1--3.2]{FZ} prove absolute continuity of the first
(mean) curvature measure with respect to surface area and obtain local
convex graph representations in $W^{2,p}$ for some $p>1$.
Corollary~\ref{cor:curvatureintro} gives singularity of the Gaussian
curvature measure and excludes a cover by $W^{2,2}$ graphs for any
full-dimensional local minimizer. The stronger hypotheses in
\cite[Theorem~1.6]{FZ} already exclude local minimality.
\end{remark}

\begin{corollary}\label{cor:energyblowup}
Suppose smooth convex bodies $K_j\to K$ in Hausdorff distance, where $K$ is a full-dimensional local minimizer of $Q$, and suppose $K_j$ contain a common ball and lie in a common bounded ball. If $II_j$ is the second fundamental form of $\partial K_j$, then
\begin{equation}\label{eq:IIblowup}
 \lim_{j\to\infty}\int_{\partial K_j}|II_j|^2\dd S=+\infty.
\end{equation}
\end{corollary}
\begin{proof}
Assume there is a subsequence, still called $K_j$, such that
$\int_{\partial K_j}|II_j|^2\dd S\le E_0<\infty$. Translate the
center of the common inner ball to the origin. Then
$\partial K_j=\{\rho_j(\omega)\omega:\omega\in\Sph^2\}$ with
$0<r\le\rho_j\le R$. On the sphere, $\rho_j^{-1}$ is the Minkowski gauge of $K_j$.
The ball bounds then give
$\|D\rho_j\|_{L^\infty}\le C(r,R)$. In a spherical chart, up to
orientation,
\[
 (II_j)_{ab}
 =\frac{\rho_j^2\sigma_{ab}+2(\rho_j)_a(\rho_j)_b
 -\rho_j(\rho_j)_{;ab}}{\sqrt{\rho_j^2+|D\rho_j|^2}}.
\]
The induced metrics and area elements are uniformly comparable to
those of the sphere. Solve for the covariant Hessian and integrate
to obtain
\[
 \|\rho_j\|_{W^{2,2}(\Sph^2)}
 \le C\bigl(1+\|II_j\|_{L^2(\partial K_j)}\bigr).
\]
Choose a subsequence with $\rho_j\rightharpoonup\rho$ weakly
in $W^{2,2}(\Sph^2)$. Hausdorff convergence also gives uniform
convergence of the radial functions. The limit $\rho$ is the radial
function of $K$.

Choose a local boundary graph of $K$. On any smaller planar patch
compactly contained in its projection, convexity and Hausdorff
convergence give graphs $x_3=f_j(x')$ of $\partial K_j$ with uniform
Lipschitz bounds and local uniform convergence to the graph of $K$.
Shrink the patch to $U$ so a common angular chart contains the
radial directions of these graphs. In its coordinates $z$, write
$X_j(z)=\rho_j(\omega(z))\omega(z)$,
$\Psi_j(z)=\pi X_j(z)$, and $t_j(z)=e_3\cdot X_j(z)$, where $\pi$
is horizontal projection. The estimates for $\rho_j$ give Lipschitz
and $W^{2,2}$ bounds for both $\Psi_j$ and $t_j$. To invert $\Psi_j$,
send $x'$ to the angular coordinates of the normalized graph point
$(x',f_j(x'))/|(x',f_j(x'))|$. Because $f_j$ has a uniform Lipschitz bound and the radius is at
least $r$, this inverse has a uniform Lipschitz bound on $U$. Applying the second-order chain rule to
$t_j=f_j\circ\Psi_j$ gives, almost everywhere,
\[
 D^2t_j=(D\Psi_j)^T(D^2f_j\circ\Psi_j)D\Psi_j
       +\sum_{k=1}^2(D_kf_j\circ\Psi_j)D^2\Psi_j^k.
\]
The bounds on $Df_j$ and $(D\Psi_j)^{-1}$ and a change of variables
give $\sup_j\|f_j\|_{W^{2,2}(U)}<\infty$. By weak compactness, the limiting graph $f$ belongs to
$W^{2,2}_{\mathrm{loc}}(U)$. The same holds in every convex graph
chart of $\partial K$. This contradicts Corollary~\ref{cor:curvatureintro}.
\end{proof}

\section{A dyadic example and an obstruction to compact replacement}\label{sec:example}
The function $H_\mathrm d$ in \eqref{eq:Hdyadicintro} joins the
values of $t^2$ at consecutive dyadic points by straight lines. The slopes $3\cdot2^n$ are positive and increase with $n$, so $H_\mathrm d$ is convex and nondecreasing. It satisfies
\begin{equation}\label{eq:dyadicfacts}
 t^2\le H_\mathrm d(t)\le\frac98t^2,\qquad
 H_\mathrm d(2t)=4H_\mathrm d(t),\qquad
 D^2H_\mathrm d=\sum_{n\in\mathbb Z}\frac32\,2^n\delta_{2^n}.
\end{equation}
The lower bound follows from the convexity of $t^2$. For the upper
bound, set $t=2^ns$ and find the maximum of $(3s-2)/s^2$ on $[1,2]$. Using $\gamma_{\mathrm{sq}}=\max\{|x_1|,|x_2|\}$ gives
\[
 \tfrac12|x'|^2\le f_\mathrm d(x')\le\tfrac98|x'|^2.
\]
On $[-1,1]^2$, the lower graph is convex and its height is at most
$H_\mathrm d(1)=1$. So $K_\mathrm d$ in \eqref{eq:Kdyadicintro} is a
compact convex body with nonempty interior. Its lower and upper faces
are separated by a positive gap. At the origin, the quadratic upper bound makes the horizontal plane
tangent to the graph. The lower bound shows that this plane exposes the
origin. For $r_j=2^{-j}$, \eqref{eq:dyadicfacts} gives
\begin{equation}\label{eq:dyadicscale}
 \frac{f_\mathrm d(r_j\xi)}{r_j^2}=f_\mathrm d(\xi)
 =H_\mathrm d(\gamma_{\mathrm{sq}}(\xi)),
\end{equation}
so Theorem~\ref{thm:main} applies at the origin.

\subsection{The explicit perturbation}
Define
\[
 \psi(t)=-16(t-1)_++26(t-2)_+-11(t-4)_++(t-8)_+.
\]
Its support lies in $[1,8]$. Its slopes are $-16$ on $(1,2)$, $10$
on $(2,4)$, and $-1$ on $(4,8)$. Integrating on these intervals gives
\[
 \int_0^\infty t\psi(t)\dd t=0,\qquad
 \int_0^\infty tH_\mathrm d'(t)\psi'(t)\dd t
 =-72+360-288=0,
\]
\begin{equation}\label{eq:explicitnorm}
 \int_0^\infty t|\psi'(t)|^2\dd t=384+600+24=1008,
 \qquad\|\psi\|_\infty=16.
\end{equation}
Take $\eta=1/64$. The slope jumps of $H_\mathrm d+\eta\psi$ remain positive, with the smallest modified jump equal to $\frac32-16\eta$. The profile is convex and agrees with the nondecreasing $H_\mathrm d$
near zero. So it is nondecreasing everywhere. By \eqref{eq:dyadicscale}, the replacement
\[
 F_j(x')=f_\mathrm d(x')+\frac{r_j^2}{64}
 \psi\left(\frac{\gamma_{\mathrm{sq}}(x')}{r_j}\right)
\]
is a rescaling of this profile composed with the square gauge. It is convex and agrees with $f_\mathrm d$ outside the annulus $r_j\le\gamma_{\mathrm{sq}}\le8r_j$. For large $j$ the annulus lies inside the base square and the replacement stays below the upper face. By \eqref{eq:explicitnorm}, the original and new vertical fibers
differ by at most $r_j^2/4$. This gives the Hausdorff bound.

For the square gauge, $\mathscr W=(-1,1)^2$ and $|D\gamma_{\mathrm{sq}}|^2=1$ almost everywhere on $\mathscr W$, so $J=4$. Equations~\eqref{eq:positiveenergy} and \eqref{eq:explicitnorm} give
\[
 \kappa=\eta^2J\int_0^\infty t|\psi'(t)|^2\dd t
 =\frac{1}{64^2}\,4\cdot1008=\frac{63}{64}.
\]
By \eqref{eq:dyadicscale}, the rescaled pair of graphs in the area calculation is exactly $P,G$. So the quadratic term is $\kappa r_j^4$, and the Taylor remainder
in \eqref{eq:areagrowth} is $O(r_j^6)$. The capacity calculation in Section~\ref{sec:mainproof} gives $o(r_j^4)$ because
$\int_{\R^2}\psi(\gamma_{\mathrm{sq}}(\xi))\dd\xi=0$. This proves the estimates in \eqref{eq:explicitrateintro}.

\subsection{Comparison with earlier criteria}
At differentiable points of the lower graph away from the origin, $H_\mathrm d$ and the active branch of $\gamma_{\mathrm{sq}}$ are locally affine, so the graph is planar there. Nonzero points on dyadic squares and on diagonals between the active
gauge branches have more than one supporting normal. The same holds at
the seams with the vertical sides. They have no second Alexandrov derivative. The upper face and the regular portions of the sides are planar as well.

The origin has a unique tangent plane. Along $e_1$, the ratio $f_\mathrm d(te_1)/t^2$ equals $1$ at $t=2^n$ and $9/8$ at $t=\frac43 2^n$, for $n\in\mathbb Z$. As $n\to-\infty$, these sequences rule out a quadratic expansion there. So Corollary~\ref{cor:alex} applies nowhere on $\partial K_\mathrm d$. The nonsmooth edges have positive length, whereas \cite[Theorem~1.6]{FZ} requires the singular boundary to have zero $\cH^1$ measure. Each regular point other than the origin lies on a planar facet whose supporting plane meets that entire facet. So the origin cannot lie in a relatively open set of regular points
with uniquely exposed supporting normals, as required in \cite[Corollary~1.8]{FZ}.

To check every possible scale at the origin, write any sequence
$s_j\downarrow0$ as $s_j=2^{-n_j}a_j$ with $1\le a_j<2$.
Choose a subsequence with $a_j\to a\in[1,2]$. Then
\eqref{eq:dyadicfacts} gives the rescaled limit
\[
 P_a(\xi)=a^{-2}H_\mathrm d(a\gamma_{\mathrm{sq}}(\xi)).
\]
Along each nonzero ray, $P_a$ is piecewise affine with positive slope jumps, so it is not two-homogeneous. The other differentiable boundary points have planar quadratic-scale tangents, which are zero after subtracting the tangent plane. The remaining points have more than one supporting normal. So no
boundary point has a positive two-homogeneous blow-up of the type
discussed in Section~\ref{sec:intro}.

\subsection{Obstruction at degenerate quadratic tangents}\label{sec:degenerate}
\begin{proposition}\label{prop:rankone}
Let $q$ be a convex quadratic polynomial on $\R^d$ whose Hessian has a nonzero null vector. If a finite convex function $G$ coincides with $q$ outside a compact set, then $G=q$ everywhere.
\end{proposition}
\begin{proof}
Choose a nonzero null vector $e$. On each line parallel to $e$, the restriction of $q$ is affine. Subtract it from the restriction of $G$ to obtain a convex function of one variable that vanishes on both tails. Its one-sided derivatives are nondecreasing and are zero on both
tails. So they are zero everywhere. The difference is zero on each
line, and these lines cover $\R^d$.
\end{proof}

For example, $q(x_1,x_2)=x_1^2$ has no nontrivial compactly supported
convex replacement. For a rank-one quadratic tangent, Proposition~\ref{prop:rankone} blocks the compact replacement used here. Corollary~\ref{cor:curvatureintro} permits rank one almost everywhere. To extend the construction, we would need a method for nonhomothetic
level sets or a nonlocal deformation that respects the affine rulings.

In $\R^n$, $n\ge3$, the same scaling gives the candidate orders $r^{n+1}$ for the area change and $r^{n+2}$ for the quadratic material-energy remainder. A descent theorem in higher dimensions would still need the moment
construction, convex gluing, and a capacity blow-up in that dimension.

\section*{Statements and Declarations}
\noindent\textbf{Funding.} Rubin Zhang acknowledges support from the China
Scholarship Council (CSC). Ning Zhang acknowledges support from the
National Natural Science Foundation of China (Grant No.~62535007) and
Huazhong University of Science and Technology (Grant No.~5003011074).\par
\noindent\textbf{Competing interests.} The authors have no relevant financial or non-financial interests to disclose.\par
\noindent\textbf{Data availability.} No external datasets were used or
generated. This article is theoretical; the explicit constructions and
calculations needed to support its conclusions are contained in the
article.\par
\noindent\textbf{Use of generative AI.} During the preparation of this
manuscript, the authors used ChatGPT (OpenAI) to assist with literature
review, manuscript drafting and language editing, organization of the
manuscript, and preliminary checks of mathematical arguments. All
AI-assisted content was reviewed and verified by the authors. The
authors take full responsibility for the mathematical results, proofs,
references, and final manuscript.\par

\end{document}